\documentclass[12pt]{amsart}
\usepackage{amsbsy,amssymb,amsmath,amsthm,amscd,amsfonts,latexsym,amstext,delarray,
amsmath,graphicx} 
\usepackage[margin=1in]{geometry}
\usepackage{color}
\input xypic

\usepackage[all]{xy}
\usepackage{setspace}

\newtheorem{thm}{Theorem}[section]
\newtheorem{prop}[thm]{Proposition}
\newtheorem{cor}[thm]{Corollary}
\newtheorem{lem}[thm]{Lemma}

\newtheorem{defn}[thm]{Definition}
\newtheorem{rem}[thm]{Remark}

\newtheorem{prob}[thm]{Question}

\numberwithin{equation}{section}

\def\Z{{\mathbb Z}}
\def\N{{\mathbb N}}
\def\R{{\mathbb R}}
\def\C{{\mathbb C}}

\def\cA{{\mathcal A}}
\def\cB{{\mathcal B}}

\def\cI{{\mathcal I}}
\def\cJ{{\mathcal J}}

\def\cM{{\mathcal M}}
\def\cN{{\mathcal N}}
\def\cO{{\mathcal O}}

\def\cW{{\mathcal W}}

\def\Aut{{\rm Aut}}

\def\fA{{\mathfrak A}}
\def\t{{\mathfrak t}}
\def\p{{\mathfrak p}}
\def\Bl{{\mathfrak B}}

\title{KMS weights on higher rank buildings}
\author{Jake Marcinek and Matilde Marcolli}
\address{Mathematics Department, Caltech, 1200 E. California Blvd. Pasadena, CA 91125, USA}
\email{jakemarcinek@gmail.com}
\email{matilde@caltech.edu}
\date{}

\begin{document}
\maketitle

\begin{abstract}
We extend some of the results of Carey--Marcolli--Rennie on modular
index invariants of Mumford curves to the case of higher rank buildings.
We discuss notions of KMS weights on buildings, that generalize the 
construction of graph weights over graph $C^*$-algebras.
\end{abstract}

\section{Introduction}

Methods of operator algebra and noncommutative geometry were applied to
Mumford curves in \cite{CMR}, \cite{ConsMa}, \cite{CLM}, \cite{CMRV}, using
graph $C^*$-algebras associated to quotients of Bruhat--Tits trees by
$p$-adic Schottky groups and boundary algebras associated to the action
of the $p$-adic Schottky group on its limit set in the conformal boundary
of the Bruhat--Tits tree. In particular, in \cite{CMR}, invariants of Mumford
curves are obtained via modular index theory on the graph $C^*$-algebra
of the quotient of the Bruhat--Tits tree by a $p$-adic Schottky group. The
modular index theory depends on the construction of KMS weights for a
suitable time evolution on the $C^*$-algebra. These are obtained via a
combinatorial equation defining {\em graph weights}. 
The goal of this paper is to develop a similar theory of KMS weights for 
higher order buildings.

\smallskip

 We focus in particular on the case of rank $2$ buildings.
In the case of buildings of type $\tilde A_2$ and their quotients by
type rotating automorphisms, a class of $C^*$-algebras generalizing the
graph $C^*$-algebras were constructed in \cite{RS1} \cite{RS2}, as higher
rank Cuntz--Krieger algebras, which generalize the usual 
Cuntz--Krieger algebras \cite{CK}. For $\Gamma$ a group of type
rotating automorphisms of an $\tilde A_2$-building $\cB$, which acts
freely on vertices with finitely many orbits, the buildings $C^*$-algebra 
of \cite{RS1} \cite{RS2} has the very natural property of being isomorphic
to the boundary algebra $C(\partial \cB)\rtimes \Gamma$ describing the
action of the group on the totally disconnected boundary at infinity
$\partial\cB$.  

\smallskip

In this paper, by considering simple
generalizations of the combinatorial equations defining graph weights,
we introduce other possible $C^*$-algebras associated to rank $2$
buildings, which generalize the Cuntz--Krieger (CK) relations of graph
$C^*$-algebras. We first recall some facts about graph weights and
we give a cohomological interpretation of the graph weight
equation. We then consider two-dimensional analogs of 
graph weights.

\smallskip

Our construction applies to an arbitrary finite CW complex (in particular
this includes the case of spherical buildings and of certain quotients
of affine buildings). The algebra we associate to $2$-dimensional
CW complex $\cB$ is just the tensor product of two graph algebras,
respectively associated to the $1$-skeleton $\cB^{(1)}$ of $\cB$ and
to a suitably defined boundary complex $\cB_\partial$, which,
respectively, account for the incidence relations in codimension
one and two. Under suitable conditions on the graphs, these are also
higher rank Cuntz--Krieger algebras, although of a 
simpler kind than those considered in \cite{RS1} \cite{RS2}.
We introduce a suitable notion of weights, the 2D CW weights,
on $2$-dimensional CW complexes that generalize the graph
weights. We construct such weights on the resulting $C^*$-algebras and we
show that they are KMS weights 
with respect to a natural time evolution. 

\smallskip

We compare the construction of KMS weights on the algebras
of $2$-dimensional CW complexes with possible constructions
of KMS weights on the higher rank CK algebras of affine $\tilde A_2$
buildings of \cite{RS1} \cite{RS2}. We present explicit
examples illustrating the general constructions. 

\smallskip

We then discuss the case of spherical buildings of rank at least three,
where a crucial result of Tits shows that such a building $\cB$ is entirely 
determined by its foundation $E_2(C)$, where $C$ is a chamber of $\cB$,
which is an amalgam of rank two buildings. This result is the key to
the classification of spherical buildings \cite{TiWe}, \cite{WeissS}. 
Higher rank affine buildings can in turn be classified in terms of their
spherical building at infinity, \cite{WeissA}. 

\smallskip

We then describe a splicing construction for graph weights and
we show that it can be applied to the 2D CW weights. We 
show that this splicing construction applied to 2D CW weights
on the generalized $m_{ij}$-gons $\Sigma_{ij}$ in the blueprint
of a higher rank spherical building can be spliced to obtain a
2D CW weight on the entire foundation $E_2(C)$.

\smallskip

The question of extending the results of \cite{CMR} from (quotients of)
Bruhat--Tits trees to higher rank buildings was posed to the second 
author by Ludmil Katzarkov, in relation to the recent work \cite{KNPS}.
While at the moment we do not see a direct connection between the
operator algebraic approach described here and the construction
of \cite{KNPS}, the present work is motivated by this longer term goal.

\smallskip

{\bf Acknowledgment.} The first author was supported by a
Summer Undergraduate Research Fellowship at Caltech.
The second author is supported by NSF grants DMS-1007207, 
DMS-1201512, PHY-1205440.

\section{Graph $C^*$-algebras, graph weights, and KMS weights}\label{graphsetup}

In this section we recall some essential aspects of the construction of
graph weights and KMS weights on graph $C^*$-algebras, as
obtained in \cite{CMR}. We also give a more geometric 
description of the combinatorial graph weight equation, in terms of
a cohomological condition.

\subsection{Graph $C^*$-algebras}
We associate to any directed
graph $E = (E^0, E^1, s, r)$ the $C^*$-algebra $C^*(E)$ generated by the 
projections $\{P_v | v \in E^0\}$ and the partial isometries $\{S_e | e \in E^1\}$, 
subject to the relations
\begin{equation}\label{graphCK1}
    S_e^* S_e = P_{r(e)}
\end{equation}
for all $e \in E^1$, and
\begin{equation}\label{graphCK2}
    P_v = \sum\limits_{s(e) = v} S_e S_e^*
\end{equation}
for every $v \in E^0$.  
We refer the reader to \cite{Kumjian} for a survey of graph $C^*$-algebras.

\smallskip

In particular, it is known \cite{KPRR}
that the graph $C^*$-algebra $C^*(E)$ of a directed finite graph $E$ with 
no sources and no sinks is a Cuntz--Krieger algebra, as defined in \cite{CK}.
These are algebras generated by partial isometries $S_a$, for a finite
alphabet $a\in \fA$, with relations 
\begin{equation}\label{CKrels}
S_a^* S_a = \sum_b A_{ab} \, S_b S_b^* , \ \ \ \  \sum_a S_a S_a^* =1 ,
\end{equation}
where $A=(A_{ab})$ is an $\#\fA \times \#\fA$-matrix with entries in $\{ 0,1\}$.

\medskip
\subsection{States, weights, and time evolutions}

We recall the notion of states and weights on $C^*$-algebras, time evolutions,
and the KMS condition for equilibrium states.

\begin{defn}\label{KMSstateDef}
A state on a unital $C^*$-algebra $\cA$ is a continuous linear functional
$\varphi: \cA \to \C$ satisfying normalization $\varphi(1)=1$ and positivity
$\varphi(x^* x)\geq 0$ for all $x\in \cA$. Let $\sigma: \R \to \Aut(\cA)$ be
a continuous $1$-parameter family of automorphisms (a time evolution).
A state $\varphi$ is a KMS$_\beta$ state, for some $\beta \in \R_+$, if
for all $x,y \in \cA$ there exists a function $F_{x,y}$ that is
analytic on the strip $\cI_\beta =\{ z\in \C \,|\, 0<\Im(z)<\beta \}$ and continuous
on the boundary $\partial \cI_\beta$, satisfying $F_{x,y}(t)=\psi(\sigma_t(x)y)$
and $F_{x,y}(t+i\beta)=\psi(y \sigma_t(x))$. 
\end{defn}

For details on the properties of KMS states, we refer the reader to the
extensive treatment in \cite{BR}. An equivalent formulation of the KMS
condition is obtained by requiring the existence of a dense subalgebra
of analytic elements, invariant under the time evolution, where the
identity $\varphi(xy)=\varphi(y \sigma_{i\beta}(x))$ holds for all $x,y$.

\medskip

Weights on $C^*$-algebras are defined as follows, see \cite{Comb}.
As in the case of states, there is a GNS representation associated to weights
on $C^*$-algebras.

\begin{defn}\label{weightdef}
A weight on a $C^*$-algebra $\cA$ is a function $\psi: \cA^+ \to [0,\infty]$,
such that $\psi(x+y)=\psi(x)+\psi(y)$ for all $x,y\in \cA^+$ and $\psi(\lambda x)=
\lambda \psi(x)$ for all $\lambda \in \R_+$ and all $x\in \cA^+$. A weight extends to
a unique linear functional $\psi: \cM_\psi \to \C$, where $\cM_\psi$ is the span
of all elements $a\in \cA^+$ with $\psi(x)<\infty$. The weight is densely defined
if $\cM_\psi$ is dense in $\cA$. The weight is lower semi-continuous if the
set $\{ x\in \cA^+ \,|\, \psi(x)\leq \lambda \}$ is closed, for all $\lambda \in \R_+$.
A non-zero weight is proper if it is both densely defined and lower semi-continuous.
\end{defn}

Let $\cN_\psi=\{ x\in \cA \,|\, \psi(x^* x)<\infty \}$, so that $\cM_\psi=\cN_\psi^* \cN_\psi$.
Suppose given a continuous $1$-parameter family $\sigma_t$
of automorphisms of $\cA$. A proper weight $\psi$ on $\cA$ is a
{\em KMS weight}, with respect to the time evolution $\sigma_t$
if $\psi$ is an equilibrium weight, $\psi\circ \sigma_t =\psi$, and, for all 
$x,y \in \cN_\psi \cap \cN_\psi^*$, there is a function $F_{x,y}$ that is
analytic on the strip $\cI_1 =\{ z\in \C \,|\, 0<\Im(z)<1 \}$ and continuous
on the boundary $\partial \cI_1$, satisfying $F_{x,y}(t)=\psi(\sigma_t(x)y)$
and $F_{x,y}(t+i)=\psi(y \sigma_t(x))$. 
Notice how this definition matches the KMS$_1$ condition for
states discussed above. 

A different way of defining the KMS condition for weights would be 
by requiring that $\psi\circ \sigma_t =\psi$ and that, 
for all $x$ in the domain of $\sigma_i$ and $xy\in \cM_\psi$,
one has $\psi(xy)=\psi(y \sigma_i(x))$.
If the weight is faithful, the time evolution $\sigma$ is uniquely determined by
$\psi$ and the KMS condition and is referred to as the 
modular group of $\psi$. 

See \cite{Comb} and \cite{Kust} for more
details on KMS weights and for the equivalence of various 
different definitions. For recent results on KMS weights on graph $C^*$-algebras,
see also \cite{Thoms}.

\medskip
\subsection{Graph weights and KMS weights on graph algebras}\label{KMSwSec}

In \cite{CMR} a construction of KMS weights on graph $C^*$-algebras
is obtained in terms of a combinatorial notion of {\em graph weights} and
the construction of explicit solutions to the corresponding graph weight
equation.

\smallskip

Let $E$ be a finite graph, $E^0 = \{v_1, \ldots, v_n\}$ with $v_{r+1}, \ldots, v_{n}$ the sinks.  As usual, let $s, r : E^1 \rightarrow E^0$ be the source and range maps.  For any vertex $v \in E^0$, we define the \emph{edge bundle} at $v$ to be the set $B_v = \{e \in E^1 | s(e) = v\}$.

\begin{defn}
    A \emph{generalized graph weight} on $E$ is a pair of $\R$-valued 
    functions $(g, \lambda)$ on $E^0$ and $E^1$, respectively, satisfying
    \begin{equation}\label{graphweight}
        g(v) = \sum\limits_{e \in B_v} \lambda(e) g(r(e))
    \end{equation}
    for each $v \in E^0$ that is not a sink.  A generalized graph weight $(g, \lambda)$ is called
    \begin{enumerate}
        \item[(i)  ]  \emph{faithful} if $g$ is never zero and
        \item[(ii)  ]  \emph{special} if $\lambda$ is constant.
    \end{enumerate}
    The word ``generalized'' is dropped and $(g, \lambda)$ is simply called a \emph{graph weight} if $g$ and $\lambda$ are nonnegative.
\end{defn}

\smallskip

Let $\mu=\mu_1 \cdots \mu_n$ denote a sequence of oriented
edges in $E^1$ with $s(\mu_{i+1})=r(\mu_i)$. The 
linear span of elements of the form $S_\mu S_\nu^*$, for
oriented paths $\mu$ and $\nu$ with $r(\mu)=r(\nu)$, is
dense in the graph $C^*$-algebra $C^*(E)$, see \cite{Kumjian}.

\smallskip

As shown in Theorem 4.5 of \cite{CMR}, 
there is a one-to-one correspondence between faithful graph weights on a
locally finite directed graph $E$ and faithful proper weights on $C^*(E)$,
with ${\rm span}\{ S_\mu S_\nu^* \}\subseteq \cM_\psi$, 
that are invariant under the gauge action defined by 
$\gamma_z(S_e)=z S_e$, for $z\in U(1)$. For the reader's convenience,
we sketch below both directions of the implication, and also the
KMS condition satisfied by these weights.

\smallskip

Let $(g,\lambda)$ be a faithful graph weight on $E$. Consider the linear functional
\begin{equation}\label{psigweight}
\psi_{(g,\lambda)} : {\rm span}\{ S_\mu S_\nu^* \} \to \C, \ \ \ 
\psi_{(g,\lambda)}(S_\mu S_\nu^*)=\delta_{\mu,\nu}\, \lambda(\nu) \, g(r(\nu)),
\end{equation}
where for a path $\nu=\nu_1\ldots\nu_n$ we set
$\lambda(\nu):= \lambda(\nu_1)\cdots \lambda(\nu_n)$. It follows from
Proposition 4.4 of \cite{CMR} that this $\psi_{(g,\lambda)}$ is a KMS
weight on $C^*(E)$, with respect to the time evolution defined on the
generators as
\begin{equation}\label{tevol}
\sigma_t(S_e)= \lambda(e)^{it}\, S_e.
\end{equation}
The KMS condition follows from the graph weight equation, the 
relations \eqref{graphCK1}, \eqref{graphCK2}, and 
$$ \sigma_t(S_\mu S_\nu^*)=\left( \frac{\lambda(\mu)}{\lambda(\nu)}\right)^{it} \, 
S_\mu S_\nu^*. $$

\smallskip

Conversely, if $\psi: {\rm span}\{ S_\mu S_\nu^* \} \to \C$ is a proper faithful
gauge invariant weight, then setting 
$$ g(v)=\psi(P_v), \ \ \ \text{ and } \ \ \ \lambda(e)=\frac{\psi(S_e S_e^*)}{\psi(S_e^* S_e)} $$
determines a faithful graph weight.

\smallskip

Thus, the question of constructing KMS weights with respect to
suitable time evolutions, is phrased in \cite{CMR} in terms of the
following combinatorial question.

\begin{prob}
    Does there exist a faithful special graph weight on $E$?  Specifically, in the case studied in \cite{CMR}, we are interested in the case $\lambda \in (0, 1)$.
\end{prob}

In \cite{CMR} a method for constructing solutions is presented, which is
adapted to the type of graphs that occur as quotients of Bruhat--Tits
trees by $p$-adic Schottky groups, namely graphs that consist of a 
finite graph (the dual graph of the special fiber of the Mumford curve)
with infinite trees attached to (some of) its vertices, \cite{Mum}. We
discuss here a cohomological method of addressing the same question.

\medskip
\subsection{Graph weights: cohomological approach}\label{constr}

The approach is as follows.  Fix $\lambda \in \R$ (or restrict to $(0,1)$ if preferred).  
Construct a chain complex and dual cochain complex, whose 0-cocycles are precisely 
the special generalized graph weights with parameter $\lambda$.  The existence of 
nontrivial special generalized graph weights is equivalent to the nontriviality of 
the 0-cohomology group, $H^0$.  We then inspect the boundary map to check 
whether such nontrivial special generalized graph weights are faithful special 
graph weights.

\smallskip

Let the 0-chains be $$C_0 = \bigoplus\limits_{v \in E^0} \R v$$ the $n$-dimensional $\R$-vector space with the vertices as a basis.  Let the 1-chains be $$C_1 = \bigoplus\limits_{v \in E^0} \R B_v$$ the $n$-dimensional $\R$-vector space with the edge bundles as a basis.  We have a chain complex
$$ 0 \longrightarrow C_1 \stackrel{\partial}{\longrightarrow} C_0 \longrightarrow 0 $$
where $$\partial(B_v) = \begin{cases} 0 &\mbox{if } B_v = \emptyset \\ \lambda \sum_{e \in B_v} r(e) - v &\mbox{otherwise.} \end{cases}$$

Now dualize to obtain the cochain complex
$$   0 \longleftarrow C^1 \stackrel{\delta}{\longleftarrow} C^0 \longleftarrow 0    $$
    
    \begin{lem}\label{0cocycles}
    The elements of the subspace $Z^0 = \ker(\delta) \leq C^0$ are the special generalized graph weights.
    \end{lem}
    
  \proof  
It is easy to see that the 0-cocycles (the subspace $Z^0 = \ker(\delta) \leq C^0$) are precisely the special generalized graph weights, because for all $v \in E^1$
\begin{equation}
    0 = \delta g(B_v) = g(\partial B_v) = \begin{cases} 0 &\mbox{if } B_v = \emptyset \\ \lambda \sum_{e \in B_v} g(r(e)) - g(v) &\mbox{otherwise} \end{cases}
\end{equation}
either gives no relation if $v$ is a sink or gives the relation in Equation \ref{graphweight} in the case where $\lambda$ is constant.
\endproof

\begin{prop}\label{H0weights}
There are nontrivial special generalized graph weights if and only if $H^0\neq 0$.
This happens if and only if $\det(\partial) = 0$. 
\end{prop}

\proof
We always have the trivial special graph weight with $g=0$.  There are nontrivial special generalized graph weights if and only if $H^0 = Z^0$ is nontrivial if and only if $H_0$ is nontrivial if and only if $\partial$ is not surjective (or equivalently not injective).  This is the case if and only if $\det(\partial) = 0$.  

With respect to the ordered bases $\{B_{v_1}, \ldots, B_{v_n}\}$ for $C_1$ and $\{v_1, \ldots, v_n\}$ for $C_0$, the matrix representation of $\partial$ is
\begin{equation}\label{special}
    M = (\lambda m_{i,j}) - (I_r \oplus 0_{n-r})
\end{equation}
where $m_{i,j}$ is the number of edges from $v_i$ to $v_j$.

Now suppose $(x_j) \in \ker(\partial)$.  Then $g \in C^0$ by $g(v_j) = x_j$ is in fact a cocycle and hence a special generalized graph weight.  Conversely, if $g \in Z^0$, then $(g(v_j)) \in \ker(\partial)$.  In particular, there exists a faithful special graph weight with parameter $\lambda$ if and only if $$\ker(M) \cap \R_{>0}^n \neq \emptyset$$ or, equivalently, 1 is an eigenvalue of $M + I$ that has an eigenvector with strictly positive components.  This second equivalence is precisely the statement of Lemma 4.6 from \cite{CMR}.
\endproof

\begin{rem}\label{notspecial}{\rm 
    This reasoning is easily generalized to (not necessarily special) graph weights.  The matrix $M$ simply becomes
    \begin{equation}
        M = \left(\sum\limits_{\substack{s(e) = v_i, \\ r(e) = v_j}} \lambda(e) \right) - (I_r \oplus 0_{n-r}).
    \end{equation} }
\end{rem}

\bigskip

\section{$2$-dimensional CW complexes}\label{2DCWsec}

In this section we consider a first approach to generalizing the previous construction from
graphs to higher dimensional combinatorial objects, with particular focus 
on the $2$-dimensional case. Here we consider
the general setting of a 2-dimensional CW complex. We introduce a generalization of
graph weights, which combine a graph weight on the 1-skeleton of the CW complex,
with a graph weight on a ``boundary graph" of the 2-dimensional complex. We then
discuss a construction of a $C^*$-algebra of the 2-dimensional complex, which combines
the graph $C^*$-algebras of the 1-skeleton and the boundary graph. 
We consider finite 2-dimensional CW complexes, oriented in the following sense.

\begin{defn}\label{2buildDef}
    Let $\cB$ be a finite 2-dimensional CW complex.  We say $\cB$ is \emph{oriented} if its 1-skeleton, $\cB^{(1)}$, is a directed graph with range and source maps $r,s : \cB^1 \rightarrow \cB^0$ and there is a map $$\nu : \cB^2 \rightarrow \bigsqcup_n \left(\prod^n \cB^1\right) / \gamma_n$$ where $\gamma_n \leq S_n$ is the subgroup of cyclic permutations, i.e. generated by $(1 2 \dots n)$.  We further require if $\nu(\sigma) = (e_1, \ldots, e_n)$, then $r(e_i) = s(e_{i+1}), i \in \Z/n\Z$.
\end{defn}

Now to any finite 2-dimensional CW complex, we associate a notion of a \emph{boundary graph}.

\begin{defn}\label{2buildDef2}
    Let $\cB = (\cB^0, \cB^1, \cB^2, r, s, \nu)$ be an oriented finite 2-dimensional CW complex.  Define the \emph{boundary graph of $\cB$} to be the directed graph $\cB_\partial$ with $\cB_\partial^0 = \cB^1$ and an edge from $e_1$ to $e_2$ for each instance that $e_2$ follows $e_1$ over all $\nu(\sigma), \sigma \in \cB^2$.  The corresponding range and source maps are denoted $\partial r, \partial s$, respectively.
\end{defn}

\medskip
\subsection{Rank 2 graph weights and 2D CW weights}

We now propose some notions of {\em 2D CW weights}, for finite 2-dimensional CW complex,
which generalize the notion of graph weights recalled above.

The idea is to consider, separately, graph weight equations for the $1$-skeleton $\cB^{(1)}$ of the 
2-dimensional CW complex
and for the boundary graph $\cB_\partial$, and then impose a relation between the edge function $\lambda$
of the graph weight of $\cB^{(1)}$ and the vertex function $\tilde\lambda$ of the graph weight for $\cB_\partial$.

\begin{defn}\label{2graphweight}
    Let $\cB$ be an oriented finite 2-dimensional CW complex.  A quadruple of nonnegative real functions 
    $(g, \tilde\lambda, \lambda ,\eta)$ on $\cB^0$, $\cB^1$, $\cB^1$ and $\cB^2$ respectively, 
    is a \emph{rank 2 graph weight} on $\cB$ if they satisfy
    \begin{align}
        g(v) &= \sum_{s(e) = v} \tilde\lambda(e) g(r(e)) \label{cond1} \\
        \lambda(e) &= \sum_{\nu(\sigma) \ni e} \eta(\sigma) \lambda(e') \label{cond2}
    \end{align}
    for all $v \in \cB^0$ and $e \in \cB^1$ where the first sum is taken over all $e \in \cB^1$ with $s(e) = v$ and the second sum is taken over all $\sigma \in \cB^2$ and all appearances of $e$ in $\nu(\sigma)$ and $e'$ is the edge following that appearance of $e$ in $\nu(\sigma)$.
\end{defn}

\begin{rem}{\rm
    Note that Condition \ref{cond2} is precisely Condition \ref{cond1} for $\cB_\partial$.  This gives us an alternate formulation of a rank 2 graph weight. }
\end{rem}

\begin{defn}\label{form2}
    A \emph{rank 2 graph weight} on the oriented finite 2-dimensional CW complex $\cB$ is a triple of nonnegative real functions $(g, \tilde\lambda, \lambda, \eta)$ on $(\cB^0, \cB^1, \cB^1, \cB^2)$ such that $(g, \tilde\lambda)$ is a graph weight on the 1-skeleton $\cB^{(1)}$ and $(\lambda, \eta)$ is a graph weight on $\cB_\partial$.
\end{defn}

In Definitions \ref{2graphweight} and \ref{form2}, we have not imposed any relation between the
solutions of Condition \ref{cond1} and \ref{cond2}. However, it is natural to require that the functions
$\lambda(e)$ and $\tilde\lambda(e)$ on $\cB^1$ are related. We consider two possible choices
of relations between these functions. 

\begin{defn}\label{Bweight0}
Let $\cB$ be an oriented finite 2-dimensional CW complex. 
A {\em tight 2D CW weight} on $\cB$ is a rank 2 graph weight, as in
Definition \ref{2graphweight} where $\lambda(e) =\tilde\lambda(e)$ for all $e\in \cB^1$. 
A {\em 2D CW weight} on $\cB$ is a rank 2 graph weight where $\lambda(e) = \tilde\lambda(e) g(r(e))$,
for all $e\in \cB^1$.
\end{defn}

We can then refer to a 2D CW weight (or tight 2D CW weight) as a triple of functions
$(g,\lambda,\eta)$ on $(\cB^0,\cB^1,\cB^2)$.
In analogy to the graph case we make the following definition.

\begin{defn}\label{Bweight1}
    A 2D CW weight (or tight 2D CW weight) $(g, \lambda, \eta)$ on an oriented finite 
    2-dimensional CW complex $\cB$ is called
    \begin{enumerate}
        \item[(i)  ]  \emph{faithful} if $g$, $\lambda$, and $\eta$ are never zero and
        \item[(ii)  ]  \emph{special} if $\eta$ is constant.
    \end{enumerate}
    If we loosen the definition of a (tight) 2D CW weight and we only require that $(g, \lambda, \eta)$ is a (possibly negative) triple of real functions, we say $(g, \lambda, \eta)$ is a \emph{generalized (tight) 2D CW weight} 
    on $\cB$.
\end{defn}

The strategy for constructing faithful special 2D CW weights (or tight 2D CW weights) is then
summarized as follows. Let $\cB$ be a finite 2-dimensional CW complex.  Starting at the top dimension and moving 
down, we obtain a similar result as Lemma 4.6 in \cite{CMR}.  We are interested in determining 
whether $\cB$ admits a faithful special (tight) 2D CW weight.  Inspired by Definition \ref{form2}, 
we first want to determine whether the graph $\cB_\partial$ admits a special graph weight.  
By \S \ref{constr}, we have a bijective correspondence between the space of faithful special 
graph weights $(\lambda, \eta)$ on $\cB_\partial$ and $\ker(M_{\cB_\partial}) \cap \R^n_+$ 
where $M_{\cB_\partial}$ is the matrix from \ref{special} corresponding to $\cB_\partial$.

\begin{rem}{\rm
    If every edge in $\cB^1$ belongs to a face in $\cB^2$, then there are no sinks in $\cB_\partial$, 
    hence $\det(M_{\cB_\partial})$ is a polynomial in $\eta$ and each generalized faithful special graph 
    weight (up to scalar multiples of $\lambda$) on $\cB_\partial$ corresponds to a root of this polynomial, 
    of which there are finitely many.  Say the unique generalized faithful graph weights are 
    $$\{(c \lambda_i, \eta_i) | 1 \leq i \leq k, c \in \R\}.$$  Since we are interested only in faithful 
    special graph weights, we do not consider any pairs with $\eta_i \leq 0$ or $c \lambda_i$ 
    not everywhere positive.}
\end{rem}

\medskip

This gives a list (finite up to positive scalar multiplication of $\lambda$) of graph weights on $\cB_\partial$.  
We may then use Remark \ref{notspecial} to check which of these assignments of $\lambda$ extends to a 
graph weight on $\cB^{(1)}$.  Here we use either the condition that the function $\lambda(e)$ itself has to
be the edge function of a graph weight on $\cB^{(1)}$ (tight 2D CW weight) or the condition that there
should be a nowhere vanishing vertex function $g(v)$ such that, $\lambda(e)/g(r(e))$ should be the
edge function of a graph weight on $\cB^{(1)}$ with vertex function $g(v)$ (2D CW weight).
The assignments that do satisfy these conditions are, respectively, the tight 2D CW weights and the
2D CW weights on $\cB$.

\medskip
\subsection{Examples of 2D CW weights and tight 2D CW weights}

\begin{figure}
\begin{center}
\includegraphics[scale=1]{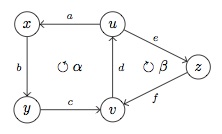}
\end{center}
\caption{The rank 2 building $\cB$.}\label{B}
\end{figure}

\begin{figure}
\begin{center}
\includegraphics[scale=1]{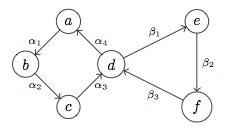}
\end{center}
\caption{The boundary graph $\cB_\partial$ of the rank 2 building $\cB$.}\label{bdryB}
\end{figure}

We now provide a simple example to illustrate the construction described above.  
Consider the oriented finite 2-dimensional CW complex 
$\cB$ given in Figure~\ref{B}.
This oriented 2-complex $\cB$ has $\cB^0 = \{u, v, x, y, z\}$, $\cB^1 = \{a, b, c, d, e\}$, and 
$\cB^2 = \{\alpha, \beta\}$ with range and source maps as shown in Figure~\ref{B}.  The two chambers in $\cB^2$ 
are attached via $(a,b,c,d)$ and $(d,e,f)$. The boundary graph $\cB_\partial$ is given in Figure~\ref{bdryB}. It
has $\cB_\partial^0 = \{a, b, c, d, e\}=\cB^1$ with range and source maps as shown in the figure.  

\begin{lem}\label{exgrwBpartial}
The faithful special graph weights on the boundary graph $\cB_\partial$ of 
Figure~\ref{bdryB}
are pairs $(\lambda,\eta)$ where $\eta$ is the unique positive root of the polynomial
$p(\eta)=1 - \eta^3 - \eta^4$ and $\lambda = C \lambda_0$, for an arbitrary $C\in \R_+^*$
and $\lambda_0$ is the function 
\begin{equation}\label{lambda0}
    \lambda_0 : (a,b,c,d,e,f) \mapsto (\eta^3, \eta^2, \eta, 1, \eta^2, \eta).
\end{equation}
\end{lem}

\proof
If $(\lambda, \eta)$ is a special graph weight on $\cB_\partial$, then 
Condition \ref{graphweight} gives the following linear system of equations
\begin{align*}
    \lambda(a) &= \eta\, \lambda(b) \\
    \lambda(b) &= \eta\, \lambda(c) \\
    \lambda(c) &= \eta\, \lambda(d) \\
    \lambda(d) &= \eta\, \lambda(a) + \eta\, \lambda(e) \\
    \lambda(e) &= \eta\, \lambda(f) \\
    \lambda(f) &= \eta\, \lambda(d)
\end{align*}
which has a nontrivial solution in $\lambda$ if and only if
\begin{equation}
    0 = \det \begin{pmatrix}
        1 & -\eta & 0 & 0 & 0 & 0 \\
        0 & 1 & -\eta & 0 & 0 & 0 \\
        0 & 0 & 1 & -\eta & 0 & 0 \\
        -\eta & 0 & 0 & 1 & -\eta & 0 \\
        0 & 0 & 0 & 0 & 1 & -\eta \\
        0 & 0 & 0 & -\eta & 0 & 1
    \end{pmatrix} = 1 - \eta^3 - \eta^4 = p(\eta).
\end{equation}
One can easily check that the polynomial $p$ has only one positive root 
$\eta_0 \in (0,1)$ so we have $\eta = \eta_0$.  Now $\lambda$ must be 
a scalar multiple of the function \eqref{lambda0}, 
say $\lambda = C \lambda_0, C \in \R_+^*$.
\endproof

When we consider the condition for 2D CW weights, we look for
pairs of functions $(g,\tilde\lambda)$ satisfying the graph weight
equation on the $1$-skeleton $\cB^{(1)}$, with the relation
$\tilde\lambda(e) g(r(e))=\lambda(e)$. We obtain the following result.

\begin{prop}\label{bwEx}
The faithful special 2D CW weights on the finite 2-dimensional CW complex of Figure~\ref{B}
are quadruples $(g,\tilde\lambda,\lambda,\eta)$ with $\eta=\eta_0$ the positive root of
$p(\eta)=1 - \eta^3 - \eta^4 = 0$, 
\begin{equation}\label{gvEx}
g : (x,y,z,u,v) \mapsto (C \eta^2, C \eta, C\eta, C\eta^{-1}, C)
\end{equation}
\begin{equation}\label{tildelambdaEx}
\lambda: (a,b,c,d,e,f) \mapsto (C\eta^3,C\eta^2,C\eta,C,C\eta^2,C\eta)
\end{equation}
and with $\tilde\lambda$ the function constant equal to $\eta=\eta_0$ on all edges.
\end{prop}

\proof By Definition \ref{Bweight0}, in order to obtain a faithful special 2D CW weight
from a graph weight $(\lambda,\eta)$ on $\cB_\partial$, we look for a faithful graph weight
$(g,\tilde\lambda)$ on $\cB^{(1)}$ with $\tilde\lambda(e) g(r(e))=\lambda(e)$. The latter condition gives
equations
\begin{align*}
\tilde\lambda(a)\, g(x) & = C \eta^3 \\
\tilde\lambda(b) \, g(y) & = C \eta^2 \\
\tilde\lambda(c) \, g(v) & = C \eta \\
\tilde\lambda(d) \, g(u) & = C \\
\tilde\lambda(e)\, g(z) & = C\eta^2 \\
\tilde\lambda(f)\, g(v) & = C \eta,
\end{align*}
while the graph weight requirement gives the equations
\begin{align*}
g(v) & = \tilde\lambda(d)\, g(u) \\
g(u) & = \tilde\lambda(a)\, g(x) + \tilde\lambda(e)\, g(z) \\
g(x) & = \tilde\lambda(b)\, g(y) \\
g(y) & = \tilde\lambda(c) \, g(v) \\
g(z) & = \tilde\lambda(f)\, g(v).
\end{align*}
These have solutions
$$ \tilde\lambda(a)=\tilde\lambda(b)=\tilde\lambda(c)=\tilde\lambda(e)=\tilde\lambda(f)=\eta, \ \ \ \tilde\lambda(d)=(\eta^3+\eta^2)^{-1} $$
$$ g(x)=C\eta^2, \ \ g(y)=C\eta, \ \ g(z)=C \eta, \ \  g(u) = C (\eta^3+\eta^2), \ \ g(v) =C. $$
Since $\eta=\eta_0$ is a root of $p(\eta)=0$, it satisfies $\eta^3+\eta^2=\eta^{-1}$, hence we
obtain the statement.
\endproof

\medskip

Similarly, we find that the solutions above fit into a larger $2$-parameter family of 
faithful 2D CW weights that have possibly different values 
$\eta(\sigma_1)\neq \eta(\sigma_2)$ for the two faces of Figure~\ref{B}.

\begin{cor}\label{Exbwnotspecial}
The faithful 2D CW weights are quadruples of functions $(g,\tilde\lambda,\lambda,\eta)$
with $\eta(\sigma_1)=\eta_1$ and $\eta(\sigma_2)=\eta_2$, where $\eta_1,\eta_2\in \R^*_+$
satisfy $\eta_1^4+\eta_2^3=1$ and with
\begin{equation}\label{gvExns}
g: (x,y,z,u,v)\mapsto (C \eta_1^2, C\eta_1, C\eta_2, C (\eta_1^3+ \eta_2^2), C)
\end{equation}
\begin{equation}\label{tildelambdaExns}
\tilde\lambda: (a,b,c,d,e,f) \mapsto (\eta_1,\eta_1,\eta_1,(\eta_1^3+\eta_2^2)^{-1},\eta_2,\eta_2)
\end{equation}
\begin{equation}\label{lambdaExns}
\tilde\lambda: (a,b,c,d,e,f) \mapsto (\eta_1^3,\eta_1^2,\eta_1,1,\eta_2^2,\eta_2).
\end{equation}
\end{cor}

\proof The argument is exactly as before with graph weight equations on $\cB_\partial$
giving 
\begin{align*}
    \lambda(a) &= \eta(\sigma_1)\, \lambda(b) \\
    \lambda(b) &= \eta(\sigma_1)\, \lambda(c) \\
    \lambda(c) &= \eta(\sigma_1)\, \lambda(d) \\
    \lambda(d) &= \eta(\sigma_1)\, \lambda(a) + \eta(\sigma_2)\, \lambda(e) \\
    \lambda(e) &= \eta(\sigma_2)\, \lambda(f) \\
    \lambda(f) &= \eta(\sigma_2)\, \lambda(d)
\end{align*}
which has a nontrivial solution in $\lambda$ if and only if
\begin{equation}\label{p12mat}
    0 = \det \begin{pmatrix}
        1 & -\eta(\sigma_1) & 0 & 0 & 0 & 0 \\
        0 & 1 & -\eta(\sigma_1) & 0 & 0 & 0 \\
        0 & 0 & 1 & -\eta(\sigma_1) & 0 & 0 \\
        -\eta(\sigma_1) & 0 & 0 & 1 & -\eta(\sigma_2) & 0 \\
        0 & 0 & 0 & 0 & 1 & -\eta(\sigma_2) \\
        0 & 0 & 0 & -\eta(\sigma_2) & 0 & 1
    \end{pmatrix} = 1 - \eta(\sigma_2)^3 - \eta(\sigma_1)^4 .
\end{equation}
The solutions are multiples $\lambda = C \lambda_0$ of the function
\begin{equation}\label{lambda012}
\lambda_0: (a,b,c,d,e,f) \mapsto (\eta_1^3, \eta_1^2, \eta_1, 1, 
\eta_2^2, \eta_2),
\end{equation}
where $\eta_1 =\eta(\sigma_1)$ and $\eta_2=\eta(\sigma_2)$,
satisfying $\eta_1^4+\eta_2^3 =1$.
We then consider the system of equations
\begin{align*}
\tilde\lambda(a) \, g(x) & = C \eta_1^3 \\
\tilde\lambda(b) \, g(y) & = C \eta_1^2 \\
\tilde\lambda(c) \, g(v) & = C \eta_1 \\
\tilde\lambda(d) \, g(u) & = C \\
\tilde\lambda(e) \, g(z) & = C\eta_2^2 \\
\tilde\lambda(f) \, g(v)  & = C \eta_2,
\end{align*}
which express the condition $\tilde\lambda(e) g(r(e))=\lambda(e)$ of
the 2D CW weights, as well as the condition
\begin{align*}
g(v) & = \tilde\lambda(d)\, g(u) \\
g(u) & = \tilde\lambda(a)\, g(x) + \tilde\lambda(e)\, g(z) \\
g(x) & = \tilde\lambda(b)\, g(y) \\
g(y) & = \tilde\lambda(c) \, g(v) \\
g(z) & = \tilde\lambda(f)\, g(v).
\end{align*}
that $(g,\tilde\lambda)$ is a graph weight on $\cB^{(1)}$.
These have solutions as in \eqref{gvExns} and \eqref{tildelambdaExns}.
This gives a $2$-paramter family of solutions depending on $C,\eta_1,\eta_2\in \R^*_+$
with the relation $\eta_1^4+\eta_2^3 =1$.
\endproof

\medskip

In the case of tight 2D CW weights, we consider solutions $(C\lambda_0,\eta_0)$
of the faithful special graph weight equation on $\cB_\partial$, as in Lemma
\ref{exgrwBpartial}, and we impose the
condition that the same function $\lambda=C\lambda_0$ extends to
a graph weight $(g,\lambda)$ on the $1$-skeleton $\cB^{(1)}$. Thus,
we can characterize the faithful special tight 2D CW weights as follows.

\begin{prop}\label{tbwEx}
The faithful special tight 2D CW weights on the finite 2-dimensional CW complex of 
Figure~\ref{B}
are of the form $(g, C \lambda_0, \eta)$, with $\eta$ the unique positive root of
$p(\eta)=1 - \eta^3 - \eta^4$, the function $\lambda_0$ as in \eqref{lambda0},
$C$ is a positive root of $q(C)=1 - C^3\eta^3 - C^4\eta^6$ and $g(v)$ is a
solution of 
\begin{align*}
    g(u) &= C \eta^3 \, g(x) + C \eta^2\, g(z) \\
    g(x) &= C \eta^2 \, g(y) \\
    g(y) &= C \eta \, g(v) \\
    g(v) &= C \, g(u) \\
    g(z) &= C \eta \, g(v) \\
\end{align*}
\end{prop}

\proof Condition \ref{graphweight} on the 1-skeleton $B^{(1)}$ gives the system of
equations above. These have a nontrivial solution in $g$ if and only if
\begin{equation}
    0 = \det \begin{pmatrix}
        1 & -C \eta^3 & 0 & 0 & -C \eta^2 \\
        0 & 1 & -C \eta^2 & 0 & 0 \\
        0 & 0 & 1 & -C \eta & 0 \\
        - C & 0 & 0 & 1 & 0 \\
        0 & 0 & 0 & - C \eta & 1
    \end{pmatrix} = 1 - C^3\eta^3 - C^4\eta^6 = q(C).
\end{equation}
Since $\eta > 0$, $q$ has a positive root.  Positive roots are the values of 
$C$ for which there is a nontrivial faithful special tight 2D CW weight
$(g, C \lambda_0, \eta_0)$.
\endproof

\medskip
\subsection{$C^*$-algebras for finite 2-dimensional CW complexes}

We consider here a class of $C^*$-algebras $C^*(\cB)$ of 2-dimensional CW complexes $\cB$,
obtained as products of graph $C^*$-algebras for the $1$-skeleton and
the boundary graph of $\cB$. 

\begin{defn}\label{CstarB0}
Let $\cB$ be an oriented finite 2-dimensional CW complex.
Let $C^*(\cB^{(1)})$ and $C^*(\cB_\partial)$ be the graph $C^*$-algebras
associated to the $1$-skeleton $\cB^{(1)}$ and the boundary graph
$\cB_\partial$. Let $C^*(\cB)=C^*(\cB^{(1)})\otimes C^*(\cB_\partial)$.
\end{defn}

In terms of generators and relations, the algebra $C^*(\cB)$
is then generated by two independent and commuting CK families,
$\{ P_v, S_e \}$ for the graph $\cB^{(1)}$ and $\{ P_e, S_{\sigma,e} \}$ for the
boundary graph $\cB_\partial$, respectively, satisfying the relations
\begin{equation}\label{CstarBrel1}
S_e^* S_e = P_{r(e)} , \ \ \ \ P_v = \sum_{e\,:\,s(e)=v} S_e S_e^* , 
\end{equation}
when $v$ is not a sink, 
\begin{equation}\label{CstarBrel2}
S_{\sigma,e}^* S_{\sigma,e}  =P_{e'}  , \ \ \ \ P_e = \sum_{\sigma\,:\,e\in \nu(\sigma)} 
S_{\sigma,e}S_{\sigma,e}^* ,
\end{equation}
where $e'$ follows $e$ in $\nu(\sigma)$.
This construction immediately suggests a natural extension to
higher ranks.

\begin{rem}\label{CK2rk}{\rm
If the finite graphs $\cB^{(1)}$ and $\cB_\partial$ have 
neither sources nor sinks, the algebras $C^*(\cB^{(1)})$ and 
$C^*(\cB_\partial)$ are Cuntz--Krieger algebras. In that case
$C^*(\cB)=C^*(\cB^{(1)})\otimes C^*(\cB_\partial)$ is a higher rank
Cuntz--Krieger algebras (in the sense of \cite{RS1}) of rank two.}
\end{rem}

\begin{defn}\label{cMB}
Let $\cM_\cB$ be the linear span of elements of $C^*(\cB)$ of the form
$S_\mu S_\nu^* S_\Omega S_\Lambda^*$, for a 
pair of multi-indices $(\mu,\nu)$ consisting 
of two paths of oriented edges in $\cB^{(1)}$ with $r(\mu)=r(\nu)$ and
a pair of multi-induces $(\Omega,\Lambda)$ consisting of two 
paths of oriented edges in $\cB_\partial$ with $r(\Omega)=r(\Lambda)$.
\end{defn}

\begin{lem}\label{tevolB}
The subspace $\cM_\cB$ is dense in $C^*(\cB)$. 
\end{lem}

\proof It is known that, for a graph algebra $C^*(E)$, the span of
the elements $S_\mu S_\nu^*$, associated to paths of oriented
edges with $r(\mu)=r(\nu)$, is dense in $C^*(E)$. In the case of 
a product of two graph algebras, we similarly have a dense
span of products $S_\mu S_\nu^* S_\Omega S_\Lambda^*$,
with $(\mu,\nu)$ and $(\Omega,\Lambda)$ respectively given
by oriented paths in the two graphs.
\endproof

\begin{lem}\label{tevolCstarB}
Let $\cB$ be an oriented finite 2-dimensional CW complex.
Suppose given functions $\eta: \cB^2 \to \R^*_+$ and $\tilde\lambda: \cB^1 \to \R^*_+$.
Setting $\sigma_t(S_e)=\tilde\lambda(e)^{it} S_e$ and $\sigma_t(S_{\sigma,e}) =\eta(\sigma)^{it} S_{\sigma,e}$
determines a time evolution on the $C^*$-algebra $C^*(\cB)$.
\end{lem}

\proof The time evolution acts on elements $S_\mu S_\nu^*$ by
$$ \sigma_t(S_\mu S_\nu^*) = \left(\frac{\tilde\lambda(\mu)}{\tilde\lambda(\nu)}\right)^{it} \, S_\mu S_\nu^*, $$
$$ \sigma_t(S_\Omega S_\Lambda^*) = \left(\frac{\eta(\Omega)}{\eta(\Lambda)}\right)^{it} \, 
S_\Omega S_\Lambda^*, $$
where $\tilde\lambda(\mu):=\tilde\lambda(e_1)\cdots \tilde\lambda(e_n)$ for an oriented path
$\mu=e_1\cdots e_n$ and $\eta(\Omega)=\eta(\sigma_1)\cdots \eta(\sigma_n)$ for an
oriented path $\Omega=(\sigma_1,e_1)\cdots (\sigma_n,e_n)$. This extends continuously to
a time evolution on the $C^*$-algebra.
\endproof

\begin{prop}\label{2DCWKMS}
Let $(g,\tilde\lambda,\lambda,\eta)$ be a rank 2 graph weight (as in Definition \ref{2graphweight})
which is faithful (the functions $g$, $\tilde\lambda$, $\lambda$, $\eta$ are nowhere vanishing).
Set
\begin{equation}\label{psi2DCW}
 \psi(S_\mu S_\nu^* S_\Omega S_\Lambda^*) =
\delta_{\mu,\nu} \delta_{\Omega,\Lambda}\, \tilde\lambda(e_1)\cdots \tilde\lambda(e_n) g(r(e_n))\,
\eta(\sigma_1)\cdots \eta(\sigma_m) \lambda(a_m'), 
\end{equation}
where $\mu=e_1 \cdots e_n$ and $\Omega =(\sigma_1,a_1)\cdots (\sigma_m,a_m)$ with
$a_m'$ following $a_m$ in $\nu(\sigma_m)$.
This uniquely defines a 
weight $\psi: \cM_\cB \to \C$ that is gauge invariant and satisfies the KMS condition with respect 
to the time evolution determined by $\sigma_t(S_e)=\tilde\lambda(e)^{it} S_e$ and 
$\sigma_t(S_{\sigma,e}) =\eta(\sigma)^{it} S_{\sigma,e}$. Conversely, given a faithful gauge invariant
weight $\psi: \cM_\cB \to \C$, with the property that the ratio $\psi(S_{\sigma,e} S_{\sigma,e}^*)/
\psi(P_e)$ only depends on $\sigma$ and not on the chosen edge $e$ in $\nu(\sigma)$, setting
 $$ g(v)=\psi(P_v), \ \ \ \tilde\lambda(e)=\frac{\psi(S_e S_e^*)}{\psi(S_e^*S_e)}, \ \ \ 
\lambda(e) = \psi(P_e), \ \ \  
\eta(\sigma)= \frac{\psi(S_{\sigma,e} S_{\sigma,e}^*)}{\psi(S_{\sigma,e}^*S_{\sigma,e})}
$$
determines a faithful rank 2 graph weight.
\end{prop}

\proof In particular we have $\psi(S_e S_e^*)=\tilde\lambda(e) g(r(e))$ and $\psi(S_{\sigma,e} S_{\sigma,e}^*)=
\eta(\sigma) \lambda(e')$, with $e'$ following $e$ in $\nu(\sigma)$.
The KMS condition implies that $\psi(S_{\sigma,e}^* S_{\sigma,e})=g(r(e))=\psi(P_{r(e)})$.
The graph weight equation $g(v)=\sum_{s(e)=v} \tilde\lambda(e) g(r(e))$ makes this compatible
with the CK relation \eqref{CstarBrel1}. Similarly for the CK relation \eqref{CstarBrel2} and the
graph weight equation $\lambda(e)=\sum_{e\in \nu(\sigma)} \eta(\sigma) \lambda(e')$. Note that
the weight \eqref{psi2DCW} is a product $\psi=\psi_1\otimes \psi_2$ of KMS weights on the
graph algebras $C^*(\cB^{(1)})$ and $C^*(\cB_\partial)$, respectively. The argument is
then analogous to the case of graph weights discussed in Proposition 4.4 and
Theorem 4.5 of \cite{CMR}.
\endproof

\medskip
\subsection{A comment on 2D CW weights and algebras}

In the construction of the algebra $C^*(\cB)$ and the KMS weights associated to 
rank 2 graph weights, there are no relations between the projectors
$P_e$ of \eqref{CstarBrel2} and the projectors in \eqref{CstarBrel1}, hence
the resulting algebra $C^*(\cB)$ is just a product of two independent CK algebras,
the graph algebras $C^*(\cB^{(1)})$ and $C^*(\cB_\partial)$. Similarly, at the level
of weights, we considered the general form of rank 2 graph weights, with no
a priori relation between the functions $\tilde\lambda$ and $\lambda$. 
It would seem more natural to require, in addition to the CK relations \eqref{CstarBrel2}
and \eqref{CstarBrel1}, that the projectors $P_e$ associated to the edges in the
boundary graph $\cB_\partial$ are related to the projections $S_e S_e^*$ in the
graph $\cB^{(1)}$. For example, a relation of the form $P_e = S_e S_e^*$ would reflect, 
at the level of KMS states, the 2D CW weight relation $\tilde\lambda(e)g(r(e))=\psi(S_e S_e^*)=
\lambda(e)=\psi(P_e)$. However, in general it is not possible to impose additional relations
on the algebra $C^*(\cB)$. In fact, doing so would correspond to taking a quotient of
$C^*(\cB)$ with respect to a two-sided ideal $\cI$ generated by the additional relations.
However, it is not always possible to have nontrivial quotients of $C^*(\cB)$.
Indeed, let $E$ be a graph without sinks, satisfying the following conditions:
\begin{enumerate}
\item every loop in $E$ has an exit 
\item given any vertex $v\in E^0$ and any infinite path $\gamma$, there is a $k\in \N$ such that
there is an oriented path from the vertex $v$ to the vertex $s(\gamma_k)$ (cofinality).
\end{enumerate}
Then it is known that graph $C^*$-algebras $C^*(E)$ is simple, see Theorem 1.23 of \cite{Tomf}.
The $C^*$-tensor product of simple $C^*$-algebras with identity is again a simple
$C^*$-algebra (Theorem 1.22.6 of \cite{Tak}). Thus, if both graphs $\cB^{(1)}$ and $\cB_\partial$
satisfy the two conditions above, the algebra $C^*(\cB)$ does not have any nontrivial two-sided ideals.
One should therefore regard the special cases of 2D CW weights and tight 2D CW weights
discussed above simply as arising from KMS weights for some special choices of time evolutions
on $C^*(\cB)$ where the phase factors $\tilde\lambda(e)$ that rotate the isometries $S_e$
are related to the values $\psi(P_e)=\lambda(e)$.

\medskip
\section{$C^*$-algebras and weights for $\tilde A_2$-buildings}

The construction described in the previous section is very simple and quite general,
and it applies to arbitrary finite $2$-dimensional CW complexes. In the
present section, we consider the case of $\tilde A_2$-buildings, for which a
different construction of a $C^*$-algebra, which is a rank two generalization
of graph algebras, exists \cite{RS1}. We discuss how one can construct
KMS weights compatible with the algebras of \cite{RS1}. We refer the
reader to \S 9 of \cite{Ron} for a general description of the affine $\tilde A_{n-1}$
buildings.

\smallskip

Let $\cB$ be a locally finite thick affine rank 2 
building of type $\tilde A_2$.  In order to move to the realm of finite 
buildings, we wish to consider quotients of $\cB$ 
by finite index $\tilde A_2$ groups.  

\smallskip

Such a building $\cB$ is a rank 2 chamber system whose \emph{chambers} are triangles.  The \emph{apartments} of $\cB$  are the subcomplexes isomorphic to the Euclidean plane tesselated by triangles.  The \emph{Weyl chambers} are the $\pi / 3$-angled sectors composed of the chambers in some apartment.  We define an equivalence relation on the sectors of $\cB$.  We say sectors $A$ and $B$ are \emph{equivalent} and write $A ~ B$ if and only if $A \cap B$ is a sector.

\smallskip

The \emph{boundary} $\Omega$ of $\cB$ is defined to be the set of equivalence classes of sectors in $\cB$.  Fix some vertex $\cO$ in $\cB$ of type $0$.  For each $\omega \in \Omega$ there is a unique sector $[\cO, \omega) \in \omega$ with vertex $\cO$, see \cite{Ron}, Theorem 9.6.  We endow $\Omega$ with the topology with the collection indexed by vertices $v$ in $\cB$
\begin{equation}
    \Omega(v) = \{\omega \in \Omega \mid v \in [\cO, \omega)\}
\end{equation}
as a base for the topology.  In this topology, $\Omega$ is a totally disconnected compact Hausdorff space.

\smallskip

Let $\Gamma$ be a group of type rotating automorphisms of $\cB$ that acts freely with finitely many orbits on $\cB^0$.  There is a natural action of $\Gamma$ on $\Omega$.  As with the graph case, where the $C^*$-algebra $C^*(E)$ is Morita equivalent to $C(\partial E) \rtimes \pi_1(E)$, we have a $C^*$-algebra $C^*(\cB / \Gamma)$ associated to $\cB / \Gamma$, which is Morita equivalent to the crossed-product algebra $C(\Omega) \rtimes \Gamma$. The latter is shown to be a higher rank Cuntz--Krieger algebra, \cite{RS1}.
We first recall the construction of this algebra, from \S 1 of \cite{RS1}.

\smallskip
\subsection{Higher rank Cuntz--Krieger algebras of $\tilde A_2$ buildings}

Consider a Coxeter complex of type $\tilde A_2$, isomorphic to the apartments in $\cB$.  Each vertex is assigned a type in $\Z/3$.  Fix a vertex of type $0$ as the origin and coordinatize the vertices by $\Z^2$ with the axes being two of the three walls meeting at $(0,0)$.  Let $\t$ be the model \emph{tile} and $\p_m$ the model parallelogram of shape $m = (m_1, m_2)$ based at $(0,0)$.  That is, $\p_m$ is the parallelogram spanned by $(0, m_2 + 1)$ and $(m_1 + 1, 0)$ and $\t = \p_{(0,0)}$.

\smallskip

Now let $\Bl_m$ be the set of type rotating isometries $\p_m \rightarrow \cB$.  Then let
\begin{enumerate}
    \item $\cW_m = \Gamma \backslash \Bl_m$,
    \item $\Bl = \cup_m \Bl_m$,
    \item $\cW = \cup_m \cW_m$.
\end{enumerate}
In the special case of $\t$, we call these sets
\begin{enumerate}
    \item $\cI = \Bl_{(0,0)}$, the set of type rotating isometries $\t \rightarrow \cB$,
    \item $A = \Gamma \backslash \cI$.
\end{enumerate}
For any shape $m \in \Z_+^2$, we have two maps $t, o : \Bl_m \rightarrow \cI$ by $t(p)(l) = p(m+l)$ and $o(p) = p\mid_\t$.

\smallskip

Next we consider two $\{0,1\}$ matrices with entries indexed by $A$.  For $a,b \in A$, let
\begin{equation}
    M_i(a,b) = \begin{cases} 1 &\mbox{if } \exists p \in \Bl_{e_i} \textrm{ such that } a = \Gamma o(p), b = \Gamma t(p) \\ 0 &\mbox{otherwise.} \end{cases}
\end{equation}

\smallskip

We then follow the construction of the $C^*$-algebra $\cA$ from \S 1 of \cite{RS1} with respect to the alphabet $A$ and transition matrices $M_1$ and $M_2$.  It is shown 
in \cite{RS1} that the words $W_m$ correspond to $\cW_m$ and the decorated words 
$\bar W_m$ correspond to $\bar \cW_m$.  Now recall the corresponding $C^*$-algebra is generated by the partial isometries $\{S_{u,v} \mid u,v \in \bar W \textrm{ and } t(u) = t(v)\}$ subject to the relations
\begin{enumerate}
    \item $S_{u,v}^* = S_{v,u},$
    \item $S_{u,v} S_{v,w} = S_{u,w},$
    \item $S_{u,v} = \sum\limits_{w\in W; \sigma(w) = e_j, \\ o(w) = t(u) = t(v)} S_{uw,vw},$ for $1 \leq j \leq r$
    \item $S_{u,u} S_{v,v} = 0$ for $u \neq v \in \bar W_0$.
\end{enumerate}

\medskip
\subsection{Weights on $\tilde A_2$ buildings and their quotients}

We now look for suitable generalizations of the graph weights equations,
similar to the general case of $2$-dimensional CW complexes considered
in the previous section, but adapted to the relations of the rank $2$
Cuntz-Krieger $C^*$-algebra $\cA$ of \cite{RS1} recalled above. 
We show that there is a very simple construction of KMS weights for
the algebra $\cA$ that closely resembles the case of graph weights
on trees.

\smallskip

The partial isometries generating $\cA$ are now parameterized by pairs 
of words corresponding to type rotating isometries of parallelograms into $\cB$.
Following the graph case, we define a state in the $C^*$-algebra first on the span of elements of the form $S T^*$ where $S$ and $T$ are partial isometries from the usual generating set, i.e. of the form $S_{u,w}$, $u, w \in W$.  Following the graph case further, we would like $\phi(S_{u,w} S_{v,x}^*) = 0$ unless $(u,w)=(v,x)$ in which case $S_{u,w} S_{v,x}^* = S_{u,u}$ by our initial projection.  We therefore first look for weights 
satisfying $\tilde g(u) = \phi(S_{u,u})$.
We now have two final projections inducing two relations on these weights:
\begin{equation}\label{tildegtildeA2}
    \tilde g(u) = \sum\limits_{\substack{w\in W; \sigma(w) = e_i, \\ o(w) = t(u)}} \tilde g(uw)
\end{equation}
for both $i=1,2$.

\smallskip

Notice that, unlike the case of the 2D CW weights considered in the
previous section, these two relations do not involve 
cells in different dimensions.  They simply show how different partial isometry 
weights relate as the embedded parallelogram is expanded to a new row of 
chambers in each of the two directions. In the higher rank case of affine $\tilde A_n$
buildings one similary expects the relations to explain how embedded 
$n$-parallelepiped weights relate as the $n$-parallelepiped is expanded 
in each of the $n$ possible directions.

\smallskip

\begin{prop}\label{KMSweightsA2tilde}
The positive cone of the $\#A$-dimensional real space 
parameterized by $\{\tilde g(v) \in \R_+ | v \in A\}$ determines solutions
of \eqref{tildegtildeA2} of the form
\begin{equation}\label{solA2tilde}
    \tilde g(u) = q^{-(m_1 + m_2)} \tilde g(o(u)),
\end{equation}
where $(m_1, m_2) = \sigma(u)$. 
\end{prop}

\proof We construct two graphs, $G_1$ and $G_2$, that have vertex set $W$ and an 
edge from $u$ to $uw$ for each $\sigma(w) = e_i$ with $o(w) = t(u)$.  Using the transition matrices $M_i$, we already know how to find solutions for graph weights.  
To search for $\tilde g$, we then can simply take the intersection of the space of 
graph weights with all edge weights equal to $1$ on $G_1$ and on $G_2$, respectively. 
As discussed in \cite{Rob}, there are exactly $q+1$ faces adjacent to any given edge.  
Thus, each vertex (word of shape $m = (m_1, m_2)$) is the source of exactly $q^2$ 
edges with distinct ranges in $G_i$, all of which are words of shape $m + e_i$.
In other words, each $G_i$ is the union of directed trees of valence $q^2$.  Since the edge weights are all set to be equal to one, we would expect the weights to decay exponentially with factor $q^{-2}$.  Graph weights on trees are very easy to construct and to match up between $G_1$ and $G_2$.  For example, one can take
$$   \tilde g(u) = q^{-(m_1 + m_2)} $$
or more generally, any function of the form \eqref{solA2tilde}.
\endproof

This case, as is clear from the fact that it is constructed using 
graph weights on trees with edge weights equal to one, corresponds to a trivial time
evolution on the algebra $\cA$. In order to see more interesting and more
general cases that correspond to non-trivial time evolutions, it is convenient to 
reformulate the construction of the algebra $\cA$ as in \S 3.4 of \cite{RS3}. 

\smallskip

We consider triangle buildings where the group $\Gamma$ acts simply transitively on the vertex
set in a type rotating way. A projective plane $P$ of order $q=p^n$, for some prime $p$, has $q^2+q+1$ points,
and $q^2+q+1$ lines $L$, with each point lying on $q+1$ lines, and each line 
containing $q+1$ points. As shown in \S 3.3 of \cite{RS3}, the incidence relations of $(P,L)$
determine a triangle presentation of a group $\Gamma=\{ a_x, \, x\in P\,|\, a_x a_y a_z=1 \}$, with the
relations occurring 
whenever the points $(x,y,z)$ satisfy $y\in \lambda(x)$, where $\lambda$ is a bijection between
the set of points and the set of lines in $P$. There is a corresponding triangle building $\cB$, 
whose vertices and edges are the Cayley graph ${\rm Cayley}(\Gamma)$ and whose chambers 
correspond to $(g, g a_x^{-1}, g a_y)$ with $g\in \Gamma$ and $(x,y,z)$ with $y\in \lambda(x)$, as above.

\smallskip

\begin{figure}
\begin{center}
\includegraphics[scale=0.6]{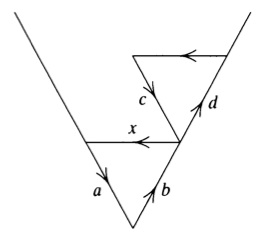}
\includegraphics[scale=0.6]{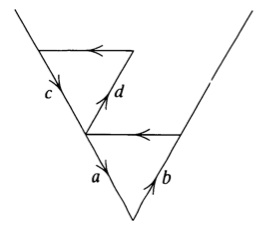}
\end{center}
\caption{Sectors defining the sets $A^\pm_{a^{-1},b}$. \label{SectorFig}}
\end{figure}

As shown in \S 3.4 of \cite{RS3}, the algebra $\cA$ can be equivalently described as generated
by two families of partial isometries $s^\pm_{a^{-1},b}$, where the pairs $(a,b)$ range over 
generators $a,b\in P$ with $b\in \lambda(a)$. There are $(q+1)(q^2+q+1)$ such elements.
Let $A^+_{a^{-1},b}$ denote the set of elements
$(c,d)$ obtained in the following way: there are $(q^2+q+1)-(q+1)=q^2$ choices of an element 
$d\notin \lambda(b)$; for each such $d$ there is then a unique $c$ satisfying $x\in \lambda(c)$
and $d\in \lambda(c)$ with $a,b,c,d,x$ in a sector as in the first diagram of 
Figure~\ref{SectorFig}.
The set $A^-_{a^{-1},b}$ is similarly defined for a sector as in the second diagram of 
Figure~\ref{SectorFig},
with $a\notin \lambda(c)$. In both cases $\# A^\pm_{a^{-1},b}=q^2$. 

Let $p_{a^{-1},b}$ denote the projection on $C(\Omega)\rtimes \Gamma$
determined by the characteristic function $\chi_{\Omega(a^{-1},b)}$ of
the clopen subset $\Omega(a^{-1},b) \subset \Omega$. 
The partial isometries $s^\pm_{a^{-1},b}$ respectively satisfy the Cuntz--Krieger
relations
\begin{equation}\label{CKplus}
s^+_{a^{-1},b}  \, {s^+}^*_{a^{-1},b} = p_{a^{-1}, b} \ \ \ \ \text{ and } \ \ \ \
{s^+}^*_{a^{-1},b} \, s^+_{a^{-1},b} = \sum s^+_{c^{-1},d}  \, {s^+}^*_{c^{-1},d},
\end{equation}
with the sum ranging over pairs $(c,d) \in A^+_{a^{-1},b}$, and
\begin{equation}\label{CKminus}
s^-_{a^{-1},b} \, {s^-}^*_{a^{-1},b} = p_{a^{-1}, b} \ \ \ \ \text{ and } \ \ \ \
{s^-}^*_{a^{-1},b} \, s^-_{a^{-1},b} = \sum s^-_{c^{-1},d} \, {s^-}^*_{c^{-1},d},
\end{equation}
summed over pairs $(c,d)\in A^-_{a^{-1},b}$. 
We also use the notation $q^\pm_{a^{-1},b} = {s^\pm }^*_{a^{-1},b} \, s^\pm_{a^{-1},b}$.

Using this presentation of the algebra $\cA$, we can reduce the construction of
KMS weights for $\cA$ to the construction of graph weights. Given a triangle building
$\cB$ constructed as above, we construct graphs $G^\pm_\cB$ with set of
vertices and edges
$$ V(G^\pm_\cB) =\{ (a,b)\in P\,:\, b\in \lambda(a) \}, $$ 
$$ E(G^\pm_\cB) = \cup_{(a,b)\in V(G^\pm_\cB)} A^\pm_{a^{-1},b}. $$
The graphs have $(q+1)(q^2+q+1)$ and $q^2$ edges out of each vertex. 

\begin{prop}\label{CKweights}
Let $\{ (g_\pm ,\lambda_\pm) \}$ be the set of faithful graph weights on the graphs $G^\pm_\cB$.
Then pairs of solutions $(g_\pm ,\lambda_\pm)$ satisfying 
\begin{equation}\label{pmmatch}
 \lambda_+(a^{-1},b) \, g_+(a^{-1},b) = \lambda_-(a^{-1},b) \, g_-(a^{-1},b),
\end{equation}
for all $(a,b)\in V(G^\pm_\cB)$, 
determine a KMS weight on the algebra $\cA$, with respect to the time evolution determined by
\begin{equation}\label{tevolpm}
 \sigma_t(s^\pm_{a^{-1},b})=\lambda_\pm(a^{-1},b)^{it} \, s^\pm_{a^{-1},b}. 
\end{equation} 
\end{prop}

\proof Set $\psi(q^\pm_{a^{-1},b})=g_\pm(a^{-1},b)$. Condition \eqref{pmmatch} ensures
that setting $$\psi(p_{a^{-1},b})=\lambda_+(a^{-1},b) \, g_+(a^{-1},b)$$ is well defined. Equivalently,
this means 
$$ \lambda_\pm (a^{-1},b) = \frac{\psi(p_{a^{-1},b})}{\psi(q^\pm_{a^{-1},b})}.  $$
By construction, $\psi$ determines a KMS weight on the CK algebras $\cA^\pm$ generated,
respectively, by the partial isometries $s^\pm_{a^{-1},b}$, with respect to the time
evolution \eqref{tevolpm}, as discussed in \S \ref{KMSwSec}. 
Condition \eqref{pmmatch} ensures that the weight $\psi$ and the time evolution  \eqref{tevolpm}
extend compatibly to the algebra $\cA$. We extend $\psi$
linearly to $\cA$ by setting 
$$ \psi(s^\alpha_\mu {s^\beta}^*_\nu) = \delta_{\mu,\nu} \delta_{\alpha,\beta} \,
\lambda_{\alpha_1}(a_1^{-1},b_1) \cdots \lambda_{\alpha_n}(a_n^{-1},b_n) \,
g_{\alpha_n}(a_n^{-1},b_n) $$
on monomials of the form $s^\alpha_\mu {s^\beta}^*_\nu$ for multi-indices 
$\mu=\mu_1,\ldots,\mu_n$, $\nu=\nu_1,\ldots,\nu_m$ with $\mu_i=(a_i,b_i)$, $\nu_j=(c_j,d_j)$,
as above, and with multi-indices $\alpha,\beta$ with $\alpha_i,\beta_j \in \{ \pm \}$. 
\endproof 

\medskip
\subsection{Triangular 2D CW weights}

Another possible construction of weights generalizing graph weights to rank $2$ buildings
of type $\tilde A_2$ can be obtained by adapting the idea of 2D CW weights discussed in
\S \ref{2DCWsec} to the triangular structure of $\tilde A_2$-buildings. 

\medskip

\begin{defn}\label{Tri2DCW}
A triangular 2D CW weight is a $5$-uple of functions $(g,\lambda,\tilde\lambda,\eta_A,\eta_B)$,
with $g$ defines on the set of vertices, $\lambda,\tilde\lambda$ on the set of edges, $\eta_A,\eta_B$
on the set of faces, satisfying
\begin{equation}\label{triangl2DCWv}
g(v) = \sum_{s(e)=v} \tilde\lambda(e) g(r(e)),
\end{equation}
\begin{equation}\label{triangl2DCW}
\lambda(e)= \sum_{\sigma\,:\,e\in\nu(\sigma)} \eta_A(\sigma) \lambda(e'') = 
\sum_{\sigma\,:\,e\in\nu(\sigma)} \eta_B(\sigma) \lambda(e') ,
\end{equation}
where $e'$ is the edge preceding $e$ in $\nu(\sigma)$ and $e''$ is the edge following $e$ in
$\nu(\sigma)$. A tight triangular 2D CW weight is as above, with $\tilde\lambda=\lambda$
and $\eta_A=\eta_B$. A triangular 2D CW weight is faithful if all the functions take strictly positive values and 
special if $\eta=\eta_A=\eta_B$ is a constant.
\end{defn}

\medskip
We show this construction in one sufficiently simple illustrative example.

\begin{prop}\label{exGamma}
The group 
\begin{equation}\label{Gammagr}
    \Gamma = \langle x_i, 0 \leq i \leq 6 \mid x_0 x_0 x_6, x_0 x_2 x_3, x_1 x_2 x_6, x_1 x_3 x_5, x_1 x_5 x_4, x_2 x_4 x_5, x_3 x_4 x_6 \rangle
\end{equation}
determines an $\tilde A_2$-building with $7$ egdes $e_i$, $7$ faces $\sigma_i$, and $1$ vertex $v$. 
Then the set of all possible special faithful tight triangular 2D CW weights on this building is a one-parameter
family given by $\{ g(v)=g\in \R^*_+, \lambda(e_i)=1/7, \eta(\sigma_i)=1/3\}$.
\end{prop}

\proof The group $\Gamma$ of \eqref{Gammagr} is an $\tilde A_2$ group of order $q=2$, hence it 
has $q^2 + q + 1 = 7$ generators.  Now consider the $\Gamma$-action on 
${\rm Cayley}(\Gamma)$.  Two triangles lie in the same $\Gamma$-orbit if and only if they have the same edge labels.  Let $\cB = {\rm Cayley}(\Gamma) / \Gamma$. The building $\cB$
has exactly $7$ edges $\{ x_0, \ldots, x_6\}$ and $7$ faces given by the relations in the 
presentation of $\Gamma$.  Equation \eqref{triangl2DCW} for a triangular 2D CW weight consists of two
sets of seven linear relations. A nontrivial solution in $\lambda$ exists if and only if
\begin{align*}
0 &= \det
\begin{pmatrix}
    1 - \eta_A  & 0 & -\eta_A & 0 & 0 & 0 & -\eta_A  \\
    0 & 1 & -\eta_A & -\eta_A & 0 & -\eta_A & 0 \\
    0 & 0 & 1 & -\eta_A & -\eta_A & 0 & -\eta_A \\
    -\eta_A & 0 & 0 & 1 & -\eta_A & -\eta_A & 0 \\
    0 & -\eta_A & 0 & 0 & 1 & -\eta_A & -\eta_A \\
    0 & -\eta_A & -\eta_A & 0 & -\eta_A & 1 & 0 \\
    -\eta_A & -\eta_A & 0 & -\eta_A & 0 & 0 & 1 \\
\end{pmatrix} \\
&= (-1 + 3 \eta_A) (1 + \eta_A + 2 \eta_A^2)^2 (-1 + 2 \eta_A^2).
\end{align*}
The matrix for $\eta_B$ is just the transpose, so one obtains the same equation.
Both give possible positive values $\{ 1/3, 1/\sqrt{2} \}$ for $\eta_A$ and $\eta_B$. By adding
rows one obtains the relations 
$$ \sum_{i=0}^6 \lambda(x_i)= 3\eta_A \sum_{i=0}^6 \lambda(x_i) = 3 \eta_B  \sum_{i=0}^6 \lambda(x_i). $$
Thus, the only case that gives rise to faithful weights is $\eta_A=\eta_B=1/3$. Any constant function
$\lambda(e_i)=\lambda >0$ is then a solution. In fact, since the matrix has rank $6$, these are the 
only solutions. For the one vertex, \eqref{triangl2DCWv} gives $g(v) = \sum_i \lambda(x_i) g(v)$,  
so this fixes the choice of $\lambda$ to be $\lambda(x_i) =\lambda= 1/7$ for all $i=0,\ldots,6$, 
while any arbitrary $g(v)=g \in \R^*_+$ will be a solution.
\endproof

These are tight 2D CW weights in the sense discussed in \S \ref{2DCWsec}, hence they correspond to
KMS weights for time evolutions on the algebra $C^*(\cB^{(1)})\otimes C^*(\cB_\partial)$ as in \S \ref{2DCWsec}.

\medskip
\section{Higher rank buildings, residues, and foundations}

We now consider cases of buildings for rank greater than two. A classification of
spherical buildings of rank at least three was given in \cite{TitsS}. A simpler proof
based on the classification of Moufang Polygons, \cite{TiWe}, is given in \cite{WeissS}. 

\smallskip

One associates to a spherical building $\cB$ an edge-colored graph $G_\cB$, whose vertex set
$V=V(G_\cB)$ is the set of chambers of $\cB$, with two chambers connected by
an edge whenever they have a common panel (codimension one faces of chambers). The set
$\cI$ of types is the set of edge coloring. If a panel has type $\cI \smallsetminus \{ i \}$, then
the corresponding edge in $E=E(G_\cB)$ is labelled with the color $i \in \cI$. Spherical buildings 
have finite apartments. Moreover, the building is thick if every panel is a face of at least three chambers. 
The rank of $\cB$ is the cardinality of $\cI$. See \cite{Ron} and \cite{WeissS} for more details.

\smallskip

A spherical building of rank two corresponds in this way to a generalized $n$-gon $G_\cB$, namely 
a connected bipartite graph with diameter $n$ and girth $2n$, where the diameter is the maximum
distance between two vertices and the girth is the length of a shortest circuit. In the thick case, vertices
of the same type have the same valence and if $n$ is odd all vertices have the same valence. Moreover,
for thick spherical buildings, $n$ is constrained to take values in the set $\{ 2,3,4,6,8 \}$, and the
valencies are also constraints, see \S 3.2 of \cite{Ron} for a detailed account. 

\medskip
\subsection{Foundations, residues, and amalgams}

More generally, a rank $N$ spherical building $\cB$ determines an $N$-partite graph $G_\cB$,
where the neighborhood of any vertex is the graph of a rank $N-1$ spherical building. Via this
reduction process, the fundamental blocks that determine the structure of rank $N$ buildings 
are identified with certain rank $2$ cases, which are special types of rank two incidence 
geometries (generalized $n$-gons):  the Moufang polygons. More precisely, 
given a subset $\cJ \subset \cI$, the $\cJ$-residue $G_{\cB,\cJ}={\rm Res}_\cJ(G_\cB)$ 
of $G_\cB$ is 
the (multi-connected) graph obtained from $G_\cB$ by removing all edges whose 
color label is not in $\cJ$. Panels correspond to $\cJ$-residues of $G_\cB$ with $\# \cJ=1$. 
The residues $G_{\cB,\cJ}$ in turn correspond to buildings $\cB_\cJ={\rm Res}_J(\cB)$. Given a chamber $C\in \cB$,
that is, a vertex $v_C\in V(G_\cB)$, the subgraph $E_2(C)\subset G_\cB$ given by the
union of the rank two residues containing $C$ is called the foundation of $\cB$. It is known that 
for thick spherical buildings of rank at least three, $\cB$ is uniquely determined by $E_2(C)$.
The foundation $E_2(C)$ is an amalgam of buildings of rank two, and can be decomposed
into a gluing of Moufang Polygons. This reduces the classification to a (difficult, but known)
classification of Moufang Polygons, obtained in \cite{TiWe}. This provides a quick sketch of
the main idea in how one obtains a classification of spherical buildings, \cite{WeissS}. 
This also suggests that, in order to construct
$C^*$-algebras, quantum statistical mechanical systems, and KMS weights, associated to the
geometry of higher rank spherical buildings, for rank at least three, it would suffice to have a
suitable construction of such objects associated to the Moufang Polygons. 

We proceed by constructing a $C^*$-algebra, obtained as described in \S \ref{2DCWsec},
associated to the foundation $E_2(C)$ of a spherical building $\cB$. We identify
$E_2(C)$ with the 2-dimensional CW-complex determined by the incidence relation of
chambers, condimension one panels and codimension two panels of $\cB$. We then
construct KMS weights on the $C^*$-algebra $C^*(E_2(C))$ obtained in this way, by
assembling tight 2D CW weights (in the sense of \S \ref{2DCWsec}) associated to the
rank two buildings given by the residues of $\cB$, whose amalgam gives $E_2(C)$.

\medskip
\subsection{Amalgams of rank two buildings}\label{rk2amalg}

Let $\cB$ be a spherical building of rank at least three. As we recalled above,
by a theorem of Tits, $\cB$ is completely determined by its foundation $E_2(C)$
which is obtained as an amalgam of rank 2 buildings, given by the union of
the residues of rank two containing the chamber $C$. 

\smallskip

A {\em blueprint} for a spherical building $\cB$ over $\cI$, with rank $N=\#\cI$, 
consists of data $\{ \Sigma_i, \Sigma_{ij} \}_{i,j\in \cI}$, where $\{ \Sigma_i \}_{i\in \cI}$ is
a {\em labeling system}, namely a system that parameterizes the $i$-residues of $\cB$.
This means that, for each residue ${\rm Res}_i(\cB)$ there is a bijection
$$ \phi_i : \Sigma_i \stackrel{\simeq}{\to} {\rm Res}_i(\cB). $$
The $\{ \Sigma_{ij} \}_{i,j\in \cI}$ are a collection of generalized $n_{ij}$-gons
with labelling by $(S_i,S_j)$,  \cite{Ron} \S 7.1.

\smallskip

In general, an amalgam of rank two buildings is given by data 
$\{ \Sigma_i, \Sigma_{ij} \}_{i,j\in \cI}$ as above such that there is 
a system of bijections 
\begin{equation}\label{phiijamalg}
 \phi_{ij}: \Sigma_i \stackrel{\simeq}{\to} {\rm Res}_i (\Sigma_{ij})
\end{equation} 
onto the $i$-th residue of $\Sigma_{ij}$, see \S 7.3 of \cite{Ron}.

\smallskip

The {\rm amalgam} $\Sigma=\amalg_{i,j} \Sigma_{ij}$ is obtained 
by gluing the generalized $n_{ij}$-gons $\Sigma_{ij}$ along the identifications
\begin{equation}\label{phiijamalg2}
 \phi_{ij}(\Sigma_i) \cong \phi_{ik} (\Sigma_i), 
\end{equation}  
that implement the $i$-adjacency relation. 
The foundation $E_2(C)$ is the amalgam of the data $\{ \Sigma_i, \Sigma_{ij} \}_{i,j\in \cI}$ 
of the blueprint.

\medskip
\subsection{Splicing graph weights}

We first discuss a construction of graph weights that reflects the operation
of splicing together two directed graphs along a common directed subgraph.

\begin{prop}\label{splicegraphs}
Let $\Gamma_1$ and $\Gamma_2$ be directed graphs, and let $\Gamma$ be a 
directed graph with embeddings $f_i :\Gamma \hookrightarrow \Gamma_i$, for
$i=1,2$, as a directed subgraph. Suppose given faithful graph weights $(g_i,\lambda_i)$
on the graphs $\Gamma_i$. Consider the graph $\Gamma_1\cup_\Gamma \Gamma_2$
obtained by gluing together the $\Gamma_i$ along the common subgraph $\Gamma$.
Then setting
\begin{equation}\label{gsplice}
g(v) =\left\{ \begin{array}{ll} g_1(v) & v \in V(\Gamma_1)\smallsetminus V(\Gamma) \\
g_2(v) & v \in V(\Gamma_2)\smallsetminus V(\Gamma) \\
g_1(v)+g_2(v)  & v\in V(\Gamma) \end{array}\right.
\end{equation}
\begin{equation}\label{lambdasplice}
\lambda(e)= \left\{ \begin{array}{ll} \lambda_1(e) & e\in E(\Gamma_1\smallsetminus \Gamma), \, r(e)\in V(\Gamma_1\smallsetminus \Gamma) \\
\lambda_2(e) & e\in E(\Gamma_2\smallsetminus \Gamma), \, r(e)\in 
V(\Gamma_2\smallsetminus \Gamma) \\[3mm]
\displaystyle{\frac{\lambda_1(e) g_1(r(e))}{g_1(r(e)) +g_2(r(e))}} &  e\in E(\Gamma_1\smallsetminus \Gamma), \, r(e)\in V(\Gamma) \\[4mm]
\displaystyle{\frac{\lambda_2(e) g_2(r(e))}{g_1(r(e)) +g_2(r(e))}} &  e\in E(\Gamma_2\smallsetminus \Gamma), \, r(e)\in V(\Gamma) \\[4mm]
\displaystyle{\frac{\lambda_1(e) g_1(r(e)) + \lambda_2(e) g_2(r(e))}{g_1(r(e)) +g_2(r(e))}} & e\in E(\Gamma)
\end{array}\right.
\end{equation}
determines a faithful graph weight $(g,\lambda)$ on the graph $\Gamma_1\cup_\Gamma \Gamma_2$.
\end{prop}

\proof Since the $(g_i,\lambda_i)$ are faithful graph weights on the $\Gamma_i$, at a vertex 
$v\in V(\Gamma)$ we have
$$ g(v)=g_1(v)+g_2(v)=\sum_{e\in E(\Gamma_1\smallsetminus \Gamma): s(e)=v} \lambda_1(e) g_1(r(e)) $$
$$  + \sum_{e\in E(\Gamma_2\smallsetminus \Gamma): s(e)=v} \lambda_2(e) g_2(r(e)) $$
$$ +  \sum_{e\in E(\Gamma): s(e)=v} (\lambda_1(e) g_1(r(e))+\lambda_2(e) g_2(r(e)). $$
The first sum, in turn, splits as two sums
$$ \sum_{e\in E(\Gamma_1\smallsetminus \Gamma): s(e)=v, r(e)\in V(\Gamma_1\smallsetminus \Gamma)} \lambda_1(e) g_1(r(e)) +
\sum_{e\in E(\Gamma_1\smallsetminus \Gamma): s(e)=v, r(e)\in V(\Gamma)} \lambda_1(e) g_1(r(e)). $$
The first of these two sums is equal to
$$ \sum_{e\in E(\Gamma_1\smallsetminus \Gamma): s(e)=v, r(e)\in V(\Gamma_1\smallsetminus \Gamma)} 
\lambda(e) g(r(e)) $$
while the second is equal to 
$$ \sum_{e\in E(\Gamma_1\smallsetminus \Gamma): s(e)=v, r(e)\in V(\Gamma)} \frac{\lambda_1(e) g_1(r(e))}
{g_1(r(e))+g_2(r(e))} \, g(r(e)) =  \sum_{e\in E(\Gamma_1\smallsetminus \Gamma): s(e)=v, r(e)\in V(\Gamma)}
\lambda(e) g(r(e)). $$
The case of the sum over $e\in E(\Gamma_2\smallsetminus \Gamma): s(e)=v$ is similar. 
The last sum is 
$$ \sum_{e\in E(\Gamma): s(e)=v} (\lambda_1(e) g_1(r(e))+\lambda_2(e) g_2(r(e)) =
\sum_{e\in E(\Gamma): s(e)=v}  \lambda(e) (g_1(r(e))+g_2(r(e)). $$
Since, if $e\in E(\Gamma)$ then also $r(e)\in V(\Gamma)$ the latter is also
$$ \sum_{e\in E(\Gamma): s(e)=v}  \lambda(e) g(r(e)). $$
Thus, the graph weight equation for the pair $(g,\lambda)$ defined as in \eqref{gsplice}, \eqref{lambdasplice}
is satisfied at all vertices $v\in V(\Gamma)$. For a vertex $v\in V(\Gamma_1\smallsetminus \Gamma)$ we
have 
$$ g(v)= g_1(v) = \sum_{e\in E(\Gamma_1\smallsetminus \Gamma), s(e)=v, r(e)\in V(\Gamma)} \lambda_1(e)
g_1(r(e)) + \sum_{e\in E(\Gamma_1\smallsetminus \Gamma), s(e)=v, r(e)\in V(\Gamma_1\smallsetminus \Gamma)} \lambda_1(e) g_1(r(e)). $$
The first sum is clearly equal to
$$ \sum_{e\in E(\Gamma_1\smallsetminus \Gamma), s(e)=v, r(e)\in V(\Gamma)} \lambda(e)
g(r(e)) $$
and the second sum is also equal to 
$$ \sum_{e\in E(\Gamma_1\smallsetminus \Gamma), s(e)=v, r(e)\in V(\Gamma_1\smallsetminus \Gamma)} \frac{\lambda_1(e) g_1(r(e))}{g_1(r(v))+g_2(r(e))} \, g(r(e)) =
\sum_{e\in E(\Gamma_1\smallsetminus \Gamma), s(e)=v, r(e)\in V(\Gamma_1\smallsetminus \Gamma)} 
\lambda(e) g(r(e)). $$
The case of a vertex in $V(\Gamma_2 \smallsetminus \Gamma)$ is analogous.
\endproof

\medskip
\subsection{Algebras and KMS weights}

Following the construction described in \S \ref{2DCWsec}, we can assign to
a foundation $E_2(C)$ of the building $\cB$ a $C^*$-algebra $C^*(E_2(C))$,
given by the higher rank Cuntz--Krieger $C^*$-algebra
$$ C^*(E_2(C)) = C^*(E_2(C)^{(1)}) \otimes C^*(E_2(C)_\partial), $$
where we identify $E_2(C)$ with a 2-dimensional CW complex and we
take $E_2(C)^{(1)}$ to be the $1$-skeleton, and $E_2(C)_\partial$ to be
the boundary complex as in \S \ref{2DCWsec}. These are, respectively,
amalgams of the $\Sigma_{ij}^{(1)}$ and the $\Sigma_{ij,\, \partial}$, with
respect to the residues $\Sigma_i$ and the identifications \eqref{phiijamalg2},
where the residues $\phi_{ij}(\Sigma_i)={\rm Res}_i (\Sigma_{ij})$ of \eqref{phiijamalg}
are seen as subcomplexes of $\Sigma_{ij}^{(1)}$. They also induce subcomplexes
of the $\Sigma_{ij,\, \partial}$.

\begin{prop}\label{spliceSigmaijKMS}
Let $\Sigma_{ij}$ be the rank two residues in the blueprint 
$\{ \Sigma_i, \Sigma_{ij} \}_{i,j\in \cI}$ of a spherical building $\cB$. 
Let $(g_{ij}, \lambda_{ij}, \eta_{ij})$ be faithful 2D CW weights constructed
on the complexes $\Sigma_{ij}$ as in \S \ref{2DCWsec}. Then the following
functions determines a faithful 2D CW weight on the foundation $E_2(C)$ of $\cB$:
\begin{equation}\label{splicegij}
g(v) = \left\{ \begin{array}{ll} g_{ij}(v) & v\in V(\Sigma_{ij}^{(1)})\smallsetminus V(\phi_{ij}(\Sigma_i)) \\
g_{ik}(v) & v\in V(\Sigma_{ik}^{(1)})\smallsetminus V(\phi_{ik}(\Sigma_i)) \\
g_{ij}(v) + g_{ik}(v) & v\in V(\phi_{ij}(\Sigma_i))=V(\phi_{ik}(\Sigma_i)) 
\end{array}\right.
\end{equation}
\begin{equation}\label{splicelambdaij}
\lambda(e) = 
\left\{ \begin{array}{ll} \lambda_{ij}(e) & e\in V(\Sigma_{ij,\,\partial})\smallsetminus E(\phi_{ij}(\Sigma_i)) \\
\lambda_{ik}(e) & e\in V(\Sigma_{ik,\,\partial})\smallsetminus E(\phi_{ik}(\Sigma_i)) \\
\lambda_{ij}(e) + \lambda_{ik}(e) & e\in E(\phi_{ij}(\Sigma_i))=E(\phi_{ik}(\Sigma_i)) ,
\end{array}\right.
\end{equation}
where $V(\Sigma_{ij,\,\partial})=E(\Sigma_{ij}^{(1)})$ and $V(\Sigma_{ik,\,\partial})=E(\Sigma_{ik}^{(1)})$, and
\begin{equation}\label{spliceetaij}
\eta(\sigma)= \left\{ \begin{array}{ll}
\eta_{ij}(\sigma) 
& (e,\sigma,e')\in E(\Sigma_{ij,\partial}\smallsetminus \phi_{ij}(\Sigma_i)_\partial), \\ &
e' \in V(\Sigma_{ij,\partial}\smallsetminus \phi_{ij}(\Sigma_i)_\partial) \\[4mm]
\displaystyle{\frac{\eta_{ij}(\sigma)\lambda_{ij}(e')}{\lambda_{ij}(e')+ \lambda_{ik}(e')}} & 
(e,\sigma,e') \in E(\Sigma_{ij,\partial} \smallsetminus \phi_{ij}(\Sigma_i)_\partial), \\ & 
e' \in V(\phi_{ij}(\Sigma_i)_\partial=\phi_{ik}(\Sigma_i)_\partial) \\[4mm]
\displaystyle{\frac{\eta_{ik}(\sigma)\lambda_{ik}(e')}{\lambda_{ij}(e')+ \lambda_{ik}(e')}} & 
(e,\sigma,e') \in E(\Sigma_{ik,\partial} \smallsetminus \phi_{ik}(\Sigma_i)_\partial), \\ & 
e' \in V(\phi_{ij}(\Sigma_i)_\partial=\phi_{ik}(\Sigma_i)_\partial) \\[4mm]
\displaystyle{ \frac{\eta_{ij}(\sigma)\lambda_{ij}(e')+\eta_{ik}(\sigma)\lambda_{ik}(e')}{\lambda_{ij}(e')+ \lambda_{ik}(e')} } & 
(e,\sigma,e') \in E(\phi_{ij}(\Sigma_i)_\partial=\phi_{ik}(\Sigma_i)_\partial),\\ &
e' \in V(\phi_{ij}(\Sigma_i)_\partial=\phi_{ik}(\Sigma_i)_\partial).
\end{array}\right. 
\end{equation}
 \end{prop}

\proof The result follows directly by applying the splicing construction of
Proposition \ref{splicegraphs} for graph weights to the pairs $(g_{ij},\tilde\lambda_{ij})$
and $(\lambda_{ij}, \eta_{ij})$, which are, respectively, faithful graph weights on the
$1$-skeleta $\Sigma_{ij}^{(1)}$  and on the boundary complexes $\Sigma_{ij,\, \partial}$,
spliced together along the residues $\Sigma_i$ via the identifications of \eqref{phiijamalg},
\eqref{phiijamalg2}. More precisely, recall that a 2D CW weight on $E_2(C)$ consists of data
$(g,\tilde\lambda,\lambda,\eta)$, where $\lambda(e)=\tilde\lambda(e) g(r(e))$ and
the pairs $(g,\tilde\lambda)$ and $(\lambda,\eta)$ are, respectively, graph weights
on the $1$-skeleton $E_2(C)^{(1)}$ and on the boundary complex $E_2(C)_\partial$. 
On the $2$-dimensional complex determined by each generalized $m_{ij}$-gon $\Sigma_{ij}$
we have a 2D CW weight, which means data $(g_{ij}, \tilde\lambda_{ij}, \lambda_{ij}, \eta_{ij})$
satisfying $\lambda_{ij}(e)=\tilde\lambda_{ij}(e) g_{ij}(r(e))$, and such that $(g_{ij},\tilde\lambda_{ij})$
is a faithful graph weight on $\Sigma_{ij}^{(1)}$ and $(\lambda_{ij}, \eta_{ij})$ is a 
faithful graph weight on $\Sigma_{ij,\, \partial}$. We apply the splicing construction to
the gluing of the $\Sigma_{ij}^{(1)}$ and of the $\Sigma_{ij,\,\partial}$ 
along the $\phi_{ij}(\Sigma_i)=\phi_{ik}(\Sigma_i)$. This gives, respectively, functions of
the form \eqref{splicegij} and 
\begin{equation}\label{splicetildelambdaij}
\tilde\lambda(e) = \left\{ \begin{array}{ll} \tilde\lambda_{ij}(e) 
& e\in E(\Sigma_{ij}^{(1)}\smallsetminus \phi_{ij}(\Sigma_i)), \\ &
r(e) \in V(\Sigma_{ij}^{(1)}\smallsetminus \phi_{ij}(\Sigma_i)) \\[4mm]
\displaystyle{\frac{\tilde\lambda_{ij}(e)g_{ij}(r(e))}{g_{ij}(r(e))+ g_{ik}(r(e))}} & 
e\in E(\Sigma_{ij}^{(1)}\smallsetminus \phi_{ij}(\Sigma_i)), \\ & 
r(e) \in V(\phi_{ij}(\Sigma_i)=\phi_{ik}(\Sigma_i)) \\[4mm]
\displaystyle{\frac{\tilde\lambda_{ik}(e)g_{ik}(r(e))}{g_{ij}(r(e))+ g_{ik}(r(e))}} & 
e\in E(\Sigma_{ik}^{(1)}\smallsetminus \phi_{ik}(\Sigma_i)), \\ & 
r(e) \in V(\phi_{ij}(\Sigma_i)=\phi_{ik}(\Sigma_i)) \\[4mm]
\displaystyle{ \frac{\tilde\lambda_{ij}(e)g_{ij}(r(e))+\tilde\lambda_{ik}(e)g_{ik}(r(e))}{g_{ij}(r(e))+ g_{ik}(r(e))} } & 
e\in E(\phi_{ij}(\Sigma_i)=\phi_{ik}(\Sigma_i)),\\ &
r(e) \in V(\phi_{ij}(\Sigma_i)=\phi_{ik}(\Sigma_i))
 \end{array}\right.
\end{equation}
and functions of the form \eqref{splicelambdaij} and \eqref{spliceetaij}. It remains to check
that the compatibility condition $\lambda(e)=\tilde\lambda(e) g(r(e))$ is satisfied, knowing
that $\lambda_{ij}(e)=\tilde\lambda_{ij}(e) g_{ij}(r(e))$ on each $\Sigma_{ij}$. 
This means checking that, given $g$ and $\tilde\lambda$ as in \eqref{splicegij} and
\eqref{splicetildelambdaij}, the function $\tilde\lambda(e) g(r(e))$ satisfies
\begin{equation}\label{checkrel}
\tilde\lambda(e) g(r(e)) = \left\{ \begin{array}{ll}
\tilde\lambda_{ij}(e) g_{ij}(r(e)) & e \in E(\Sigma_{ij}^{(1)}\smallsetminus \phi_{ij}(\Sigma_i)) \\
\tilde\lambda_{ik}(e) g_{ik}(r(e)) & e \in E(\Sigma_{ik}^{(1)}\smallsetminus \phi_{ik}(\Sigma_i)) \\
\tilde\lambda_{ij}(e) g_{ij}(r(e)) + \tilde\lambda_{ik}(e) g_{ik}(r(e)) &
e \in E(\phi_{ij}(\Sigma_i))=E(\phi_{ik}(\Sigma_i))
\end{array}\right. 
\end{equation}
This is indeed the case, by direct inspection, combining  \eqref{splicelambdaij} and \eqref{spliceetaij}.
\endproof

\smallskip

\begin{rem}{\rm
The splicing construction of Proposition \ref{spliceSigmaijKMS}
works for 2D CW weights, but would not work for tight 2D CW weight,
because the different splicing conditions \eqref{splicelambdaij} and \eqref{splicetildelambdaij}
would not allow for the tight matching condition $\lambda(e)=\tilde\lambda(e)$. }
\end{rem}

\end{document}